\documentclass[a4paper,10pt]{article}
\usepackage{paper-en}
\usepackage{hyperref}

\def\thetitle{Gromov's tori are optimal}
\def\theauthors{Anton Petrunin}

\hypersetup{colorlinks=true,
citecolor=black,
linkcolor=black,
anchorcolor=black,
filecolor=black,
menucolor=black,
urlcolor=black,
pdftitle={\thetitle},
pdfauthor={\theauthors}
}

\begin{document}

\title{\thetitle}
\author{\theauthors}
\date{}
\maketitle

\begin{abstract}
We give an optimal bound on normal curvatures of immersed $n$-torus in a Euclidean ball of large dimension. 
\end{abstract}

\section{Introduction}

Let us denote by $\BB^q$ the closed unit ball in $\RR^q$ centered at the origin.
Further, $\TT^n$ will denote the $n$-dimensional torus --- the smooth manifold diffeomorphic to the product of $n$ circles.

This note is inspired by examples of embeddings $\TT^n\hookrightarrow\BB^q$ for large $q$ with constant normal curvatures $K_n=\sqrt{3\cdot n/(n+2)}$.
In other words, any geodesic in the torus has constant curvature $K_n$ as a curve in $\RR^q$.
These examples were found by Michael Gromov among geodesic subtori in Clifford's tori 
\cite[2.A]{gromov3}, \cite[1.1.A.]{gromov2}.
In particular, Gromov's tori have flat induced metrics.
(Recall that Clifford's torus is a product of $m$  circles of radius ${1}/{\sqrt{m}}$ in $\RR^{2\cdot m}$;
its normal curvatures lie in the range $[\,1,\sqrt{m}\,]$.)

Gromov's examples lead to the following surprising facts:
\textit{any closed smooth manifold $L$ admits a smooth embedding into $\BB^q$ for large $q$ with normal curvatures less than $\sqrt{3}$;}
\textit{moreover, the induced Riemannian metric on $L$ can be chosen to be proportional to any given metric $g$;}
see \cite[1.D]{gromov3} and \cite[1.1.C]{gromov2}.

The next theorem implies that Gromov's tori have the best upper bound on normal curvatures;
in particular, the $\sqrt3$-bound is optimal.

\begin{thm}{Theorem}
\label{thm:main}
Suppose $\TT^n$ is smoothly immersed in $\BB^q$.
Then its maximal normal curvature is at least 
\[\sqrt{3\cdot \tfrac{n}{n+2}}.\]
\end{thm}

To make the statement more exact, we need one more notation.
Assume that $L$ is
a smooth $n$-dimensional manifold immersed in $\RR^q$;
we will always assume that $L$ is equipped with the induced Riemannian metric.
Let us denote by $\T_x$ and $\N_x$ the tangent and normal spaces of $L$ at~$x$.

Recall that the \emph{second fundamental form} $\II$ at $x$ is a symmetric quadratic form on $\T_x$ with values in $\N_x$.
It is uniquely defined by the identity $\II(\vec v,\vec v)\z\equiv \gamma''_{\vec v}(0)$,
where $\vec v\in \T_x$ and $\gamma_{\vec v}$ is an $L$-geodesic that starts at $x$ with initial velocity vector~$\vec v$.

Given $x\in L$,
denote by $\zh(x)$ the average value of $|\II(\vec u,\vec u)|^2$ for $\vec u\in\T_x$ such that $|\vec u|=1$.
Since $K(\vec u)=|\II(\vec u,\vec u)|$ is the normal curvature in the direction $\vec u$,
we have that
$\zh(x)$ is the average value $K^2(\vec u)$.
(The Cyrillic zhe $\zh$ is used since it resembles $K^2$.)

\begin{thm}{Theorem}
Suppose $\TT^n$ is smoothly immersed in $\BB^q$.
Let us equip $\TT^n$ with the induced Riemannian metric;
so we can take average values with respect to the induced volume.

\begin{subthm}{strong:2D}
If $n=2$, then the average value of $\zh$ is at least $\tfrac32$.
\end{subthm}

\begin{subthm}{strong:flat}
If the metric on $\TT^n$ is flat, then the average value of $\zh$ is at least $3\cdot \tfrac{n}{n+2}$.
\end{subthm}

\begin{subthm}{strong:sphere}
If the image of $\TT^n$ lies in $\partial\BB^q$, then  $\zh\ge 3\cdot \tfrac{n}{n+2}$ at some point of $\TT^n$.
\end{subthm}

\begin{subthm}{strong:4D}
If $n\le 4$, then $\zh\ge 3\cdot \tfrac{n}{n+2}$ at some point of $\TT^n$.
\end{subthm}

\begin{subthm}{strong:main}
If the normal curvatures of $\TT^n$ do not exceed $2$, then $\zh\ge 3\cdot \tfrac{n}{n+2}$ at some point of $\TT^n$.
\end{subthm}

\end{thm}

Note that part \ref{SHORT.strong:main} is a stronger version of \ref{thm:main}.
The remaining statements \ref{SHORT.strong:2D}--\ref{SHORT.strong:4D} are stronger versions of \ref{thm:main} in some partial cases.
All this follows since the normal curvature in some direction at $x$ is at least $\sqrt{\zh(x)}$.

All proofs  use our version of the Gauss formula; see below.
The proofs of \ref{SHORT.strong:sphere}--\ref{SHORT.strong:main} use in addition that 
\textit{the torus does not admit a metric with positive scalar curvature} \cite[Corollary A]{gromov-lawson}.

\begin{thm}{Open question}
Is it true that for any smooth immersion $\TT^n\looparrowright\BB^q$ the inequality $\zh\ge3\cdot \tfrac{n}{n+2}$ holds at some point? 
\end{thm}

\begin{thm}{Open question}
Suppose $\RP^n\looparrowright\BB^q$ is a smooth immersion.
Is it true that its normal curvature is at least $\sqrt{2\cdot n/(n+1)}$ in  some direction? 
\end{thm}

The last question asks if \textit{the Veronese embedding is optimal}.
Recall that the Veronese embedding $\RP^n\hookrightarrow\BB^q$ has all normal curvatures $\sqrt{2\cdot n/(n+1)}$; here $q= \tfrac12\cdot(n+1)\cdot(n+2)$.
An analogous question for immersions into unit spheres is open as well \cite{445819}.
You may also ask for the optimal constant for your favorite closed manifold.

Gromov's construction gives many smooth immersions $\TT^n\looparrowright\BB^q$ with normal curvatures $\sqrt{3\cdot n/(n+2)}$ in all directions.
The proof of our theorem implies that such immersions send torus to the boundary of $\BB^q$ and it tell exactly what is its second fundamenal form.
Nevertheless the following question is open.

\begin{thm}{Open question}
Can one describe all smooth immersions $\TT^n\looparrowright\BB^q$ with normal curvatures $\sqrt{3\cdot n/(n+2)}$ in all directions?
\end{thm}

\section{Gauss formula}

Recall that $L$ is a smooth $n$-dimensional manifold immersed in $\RR^q$.
Given $p\in L$,
denote by $\Sc(p)$ and $H(p)$
the scalar curvature and the mean curvature vector at $p$.

The following version of the Gauss formula plays a central role in all proofs;
it is used instead of the formula in \cite[5.B]{gromov1}.

\begin{thm}{Gauss formula}\label{formula:gauss}
The following identity
\[\Sc=\tfrac32\cdot |H|^2-\tfrac{n\cdot (n+2)}{2}\cdot\zh\]
holds for any smooth $n$-dimensional immersed manifold in a Euclidean space.
\end{thm}

\parit{Proof.}
Choose a point $p\in L$.

Assume $\codim L=1$.
Denote by $k_1,\dots,k_n$ the principal curvatures of $L$ at $p$.
Note that
\[|H|^2= \sum_ik_i^2+2\cdot\sum_{i<j}k_i\cdot k_j.\]
Further, 
\[
n\cdot (n+2)\cdot\zh
=
3\cdot \sum_i k_i^2+2\cdot \sum_{i<j} k_i\cdot k_j.
\]
The last identity follows since $\zh$ is the average value of $\left(\sum_i k_i\cdot x_i^2\right)^2$ on the unit sphere $\SSS^{n-1}\subset\RR^n=\T_p$;
here $(x_1,\dots,x_n)$ are the standard coordinates in~$\RR^n$.
One has to take into account that the following functions have unit average values:
$\tfrac13\cdot n\cdot (n+2)\cdot x_i^4$ and $n\cdot (n+2)\cdot x_i^2\cdot x_j^2$ for $i\ne j$.

By the standard Gauss formula,
\[\Sc=2\cdot\sum_{i<j}k_i\cdot k_j.\]
It remains to rewrite the right-hand side using the expressions for $|H|^2$ and $\zh$.

If $\codim L =k>1$, then the second fundamental form at $p$ can be presented as a direct sum of $k$ real-valued quadratic forms $\II_1\oplus\dots\oplus \II_k$;
that is,
\[\II=e_1\cdot\II_1+\dots+e_k\cdot \II_k,\]
where $e_1,\dots, e_k$ is an orthonormal basis of $\N_p$.
Denote by $\Sc_i$, $H_i$, and $\zh_i$ the values associated with $\II_i$.
From above, we get
\[\Sc_i=\tfrac32\cdot |H_i|^2-\tfrac{n\cdot (n+2)}{2}\cdot\zh_i\]
for each $i$.

Note that 
\[
\Sc=\sum_i\Sc_i,
\qquad
|H^2|=\sum_i|H_i|^2,
\quad\text{and}\quad
\zh=\sum_i\zh_i.
\]
Hence the general case follows.
\qeds

\parit{Remark.}
A more direct proof of this formula can be obtained using the so-called \emph{extrinsic curvature tensor} which is defined by
$\Phi(\vec x,\vec y,\vec v,\vec w)\df\langle\II(\vec x,\vec y),\II(\vec v,\vec w)\rangle$;
the necessary properties of this tensor are discussed in \cite{petrunin}.
As a bonus, one gets an explicit expression for the second fundamental forms of all Gromov's tori.

\section{Special cases}

The following statement appears in the book of Yuri Burago and Viktor Zalgaller \cite[Theorem~28.2.5]{burago-zalgaller};
it generalizes the result of István Fáry about average curvature of a curve in the unit ball \cite{fary,tabachnikov}, but the proof is essentially the same.

\begin{thm}{Lemma}\label{lem:av(H)}
Let $L$ be a closed $n$-dimensional manifold that is smoothly immersed in $\BB^q$.
Then the average value of $|H|$ on $L$ is at least $n$.
\end{thm}

\parit{Proof.}
Consider the function $u\:x\mapsto \tfrac12\cdot |x|^2$ on $L$.
Note that 
\[(\Delta u)(x)=n+ \langle H(x),x\rangle.\]
It follows that the average value of $\langle H(x),x\rangle$ is $-n$.
Since $|x|\le1$, we get the result.
\qeds

\parit{Proof of \ref{strong:2D}.}
By \ref{lem:av(H)}, the average value of $|H|^2$ is at least 4.
Further, by the Gauss--Bonnet formula, the scalar curvature (which is twice the Gauss curvature in this case) has zero average.
Therefore \ref{formula:gauss} implies the statement.
\qeds

\parit{Proof of \ref{strong:flat}.}
By \ref{lem:av(H)},
the average value of $|H|^2$ is at least $n^2$.
Since $\Sc\equiv0$, it remains to apply \ref{formula:gauss}.
\qeds

\parit{Proof of \ref{strong:sphere}.} 
Since the image lies in the unit sphere, we have that $|H|^2$ is at least $n^2$ at each point.
Since $\TT^n$ does not admit a metric with positive scalar curvature \cite[Corollary A]{gromov-lawson}, we have $\Sc(x)\le 0$ at some point $x$.
It remains to apply \ref{formula:gauss} at~$x$.
\qeds

\section{Main case}

The following lemma is an easy corollary of the bow lemma of Axel Schur \cite{shur,petrunin-zamora}.
It explains how we use the assumption on normal curvatures in \ref{strong:main}.
If $n\le4$, then the proof works without this assumption.

\begin{thm}{Lemma}\label{lem:trivial}
Let $L$ be a complete smooth manifold smoothly immersed in $\BB^q$.
Suppose its normal curvatures are at most $2$.
Given $x\in L$, denote by $\beta=\beta(x)$ the angle between vector $x$ and the normal space $\N_x$.
Then $|x|\le \cos\beta$.
\end{thm}

\begin{wrapfigure}{r}{39 mm}
\vskip-3mm
\centering
\includegraphics{mppics/pic-10}
\vskip0mm
\end{wrapfigure}

\parit{Proof.}
Let $\xi$ be a tangent direction at $x$ such that $\measuredangle(x,\xi)\z=\tfrac\pi2-\beta$.
In the plane spanned by $x$ and~$\xi$, choose a unit-speed circle arc $\sigma$ from $0$ to $x$ that comes to $x$ in the direction opposite to $\xi$;
extend $\sigma$ after $x$ by a unit-speed semicircle $\tilde\gamma$ with curvature $2$ in such a way that the concatenation $\sigma*\tilde\gamma$ is an arc of a $C^1$-smooth convex plane curve; see the figure.

Observe that if $|x|> \cos\beta$, then $\tilde\gamma$ leaves $\BB^q$; that is, $|\tilde\gamma(t_0)|>1$ for some $t_0$.

Let $\gamma$ be the unit-speed geodesic in $L$ that runs from $x$ in the direction~$\xi$.
Note that curvatures of $\sigma*\gamma$ do not exceed the curvatures of $\sigma*\tilde\gamma$ at the corresponding points.
Applying the bow lemma for $\sigma*\gamma$ and $\sigma*\tilde\gamma$, we get $|\gamma(t_0)|\ge |\tilde\gamma(t_0)|$.
It follows that $L$ does not lie in $\BB^q$ --- a contradiction.
\qeds

Let $g$ be a Riemannian metric on $\TT^n$.
Suppose $n\ge 3$, and $u\:\TT^n\to \RR$ is a smooth positive function.
Recall that
\[\left(\Sc\cdot u-4{\cdot}\tfrac{n-1}{n-2}{\cdot}\Delta u\right)\cdot u^{-\frac{n+2}{n-2}}\]
is the scalar curvature of the metric $u^{\frac{4}{n-2}}\cdot g$;
see for example \cite[6.3]{aubin}.
Since \emph{any Riemannian metric on $\TT^n$ has nonpositive scalar curvature at some point} \cite[Corollary A]{gromov-lawson}, we get the following.

\begin{thm}{Claim}\label{clm:sc-lap}
For any Riemannian metric on $\TT^n$
and any positive smooth function $u\:\TT^n\to\RR$, the function 
\[\Sc\cdot u-4{\cdot}\tfrac{n-1}{n-2}{\cdot}\Delta u\]
returns a nonpositive value at some point.
\end{thm}

\parit{Proof of \ref{strong:4D} and \ref{strong:main}.}
The case $n=2$ follows from \ref{strong:2D};
so we can assume that $n\ge 3$.
Consider the function $u\:x\mapsto \exp(-\tfrac k2\cdot|x|^2)$ for $x\in \TT^n$.

We will apply the following formula
\[\Delta(\phi\circ f)=(\phi'\circ f)\cdot \Delta f+(\phi''\circ f)\cdot|\nabla f|^2\]
to $f\:x\mapsto \tfrac12\cdot |x|^2$ and $\phi\:y\mapsto \exp(-k\cdot y)$; so $u=\phi\circ f$.

Set $r(x)=|x|$, $\alpha(x)=\measuredangle (H(x),x)$, and $\beta(x)$ as in \ref{lem:trivial}.
Note that 
\[\beta\le \alpha\le \pi-\beta.
\eqlbl{eq:beta=<alpha}\]
Observe that
\[\Delta f=|H|\cdot r\cdot \cos\alpha+n,
\quad
|\nabla f|=r\cdot \sin\beta,
\quad
\phi'=-k\cdot\phi,
\quad
\phi''=k^2\cdot \phi.
\]
Therefore
\[
\begin{aligned}
\Delta u
&=
u\cdot[-k\cdot|H|\cdot r\cdot \cos\alpha
-k\cdot n
+k^2\cdot r^2\cdot \sin^2\beta].
\end{aligned}
\]

Recall that
$\Sc=-\tfrac{n\cdot (n+2)}{2}\cdot\zh+\tfrac32\cdot |H|^2$; see \ref{formula:gauss}.
By \ref{clm:sc-lap}, the function
\begin{align*}
\Sc\cdot u-4{\cdot}\tfrac{n-1}{n-2}\cdot \Delta u
&=u\cdot\biggl[
-\tfrac{n\cdot (n+2)}{2}\cdot\zh
+\tfrac32\cdot |H|^2
+4{\cdot}\tfrac{n-1}{n-2}\cdot k\cdot|H|\cdot r\cdot \cos\alpha+
\\
&\qquad\qquad\qquad\qquad\qquad\qquad+4{\cdot}\tfrac{n-1}{n-2}\cdot(k\cdot n-k^2\cdot r^2\cdot \sin^2\beta)\biggr]
\end{align*}
returns a nonpositive value at some point $x\in \TT^n$.

Choose 
\[k=\tfrac34\cdot\tfrac {n-2}{n-1}\cdot n,
\quad\text{so}\quad
n=\tfrac43\cdot\tfrac {n-1}{n-2}\cdot k.\]
At the point $x$, we have
\[
\begin{aligned}
\tfrac{n\cdot (n+2)}{2}\cdot\zh
&\ge \tfrac32\cdot(|H|+n\cdot r\cdot \cos\alpha)^2
-\tfrac32\cdot n^2\cdot r^2\cdot \cos^2\alpha+
\\
&\qquad\qquad\qquad\quad+3\cdot n^2
-\tfrac94\cdot\tfrac{n-2}{n-1}\cdot n^2\cdot r^2\cdot \sin^2\beta\ge\tfrac32\cdot n^2.
\end{aligned}
\eqlbl{eq:main}
\]

Indeed,
by \ref{eq:beta=<alpha},
$\cos^2\alpha+\sin^2\beta\z\le 1$.
If $n\le 4$, then  $\tfrac32\ge \tfrac94\cdot\tfrac{n-2}{n-1}$;
therefore the last inequality follows, and it proves \ref{strong:4D}.

Further, if $n\ge 5$, then for the last inequality in \ref{eq:main}
we need to use in addition that $r^2+\sin^2\beta\le 1$ which follows from \ref{lem:trivial}.
Hence \ref{strong:main} follows.
\qeds

\parbf{Acknowledgments.}
I want to thank Michael Gromov and Nina Lebedeva for help and encouragement.
This work was done at Euler International Mathematical Institute of PDMI RAS;
it was partially supported by the National Science Foundation, grant DMS-2005279,
and the Ministry of Education and Science of the Russian Federation, grant 075-15-2022-289.

{\sloppy
\printbibliography[heading=bibintoc]
\fussy
}

\end{document}